\documentclass{amsart}
\usepackage{amsthm,amsfonts,amsmath,amscd,amssymb,latexsym,epsfig}
\newcommand{\gl}{{\mathfrak g \mathfrak l}}

\renewcommand{\u}{{\mathfrak u}}

\newcommand{\cx}{{\mathbb C}}
\newcommand{\diag}{\operatorname{diag}}

\newcommand{\tr}{\operatorname{tr}}
\newcommand{\Res}{\operatorname{Res}}
\newcommand{\re}{\operatorname{Re}}
\newcommand{\im}{\operatorname{Im}}

\newcommand{\Lie}{\operatorname{Lie}}

\newcommand{\Ker}{\operatorname{Ker}}

\newcommand{\Jac}{\operatorname{Jac }}
\newcommand{\sgn}{\operatorname{sgn }}
\newcommand{\Pic}{\operatorname{Pic }}

\numberwithin{equation}{section}

\newtheorem{theorem}{Theorem}[section]

\newtheorem{lemma}[theorem]{Lemma}

\newtheorem{corollary}[theorem]{Corollary}

\newtheorem{proposition}[theorem]{Proposition}

\theoremstyle{remark}

\newtheorem{remark}{Remark}[section]

\newtheorem{definition}{Definition}[section]

\newtheorem{example}{Example}[section]

\newcommand{\oC}{{\mathbb{C}}}

\newcommand{\oP}{{\mathbb{P}}}

\newcommand{\oR}{{\mathbb{R}}}

\newcommand{\oZ}{{\mathbb{Z}}}

\newcommand{\sA}{{\mathcal{A}}}   % connections
\newcommand{\sB}{{\mathcal{B}}}
\newcommand{\sG}{{\mathcal{G}}}   % gauge transformations

\newcommand{\sM}{{\mathcal{M}}}   % moduli space

\newcommand{\sO}{{\mathcal{O}}}

\newcommand{\sR}{{\mathcal{R}}}

\newcommand{\sV}{{\mathcal{V}}}

\newcommand{\fG}{{\mathfrak{g}}}

\newcommand{\fU}{{\mathfrak{u}}}

\begin{document}

\title[Reducible spectral curves]{Reducible spectral curves and the hyperk\"ahler geometry of adjoint orbits}
\author{Roger Bielawski}
\address{School of Mathematics\\
University of Edinburgh\\
Edinburgh EH9 3JZ\\ Scotland.}
%\date{\today}

%\thanks{}
%\dedicatory{}

\subjclass[2000]{53C26; 14L30; 14H42; 14H70}

\begin{abstract} We study the hyperk\"ahler geometry of a regular semisimple adjoint orbit of $SL(k,\cx)$ via the algebraic geometry of the corresponding reducible spectral curve. \end{abstract}

\maketitle

\thispagestyle{empty}

It is by now well-known that adjoint orbits of complex semisimple Lie groups admit hyperk\"ahler structures. Among several constructions of such structures, it is the one given by Kronheimer \cite{Kr}, later extended by Biquard \cite{Biq} and by Kovalev \cite{Kov}, that stands out. Kronheimer's hyperk\"ahler structures are algebraic and the metrics are complete or very close to being complete. They have found quite a few applications in representation theory. Yet, despite their nice properties, they remain quite mysterious. Not only are they not known explicitly, but there is also no good description of the relation between the hyperk\"ahler geometry of a semisimple orbit of $G^\cx$ and the K\"ahler geometry of the corresponding orbit of the compact group $G$. There are some exceptions: notably, Hermitian symmetric spaces studied in detail by Biquard and Gauduchon \cite{BG}.
\par
Kronheimer's construction is based on Nahm's equations: a system of Lie algebra valued ODE's used first to describe  moduli spaces of magnetic monopoles. It is well known that the Nahm equations for $\u(k)$ correspond to a linear flow on the Jacobian of an algebraic curve $S$. In Kronheimer's construction for regular semisimple orbits $SL(k,\cx)$, the underlying curve $S$ is a union of rational curves, without multiple components and with every pair of components intersecting in a pair of points. In this paper, we study  the hyperk\"ahler geometry of these orbits via the algebraic geometry of such reducible curves $S$. To make the algebraic geometry as simple and as explicit as possible, we restrict ourselves to the (generic) case of $S$ having only ordinary double points as singularities. For such an $S$, we define and give a (polynomial) formula for the theta function, analogous to one given by Mumford \cite{mum} for singular curves with only one rational component. We also give a formula for the inverse of the Abel mapping as well as one for  Hitchin's metric on spaces of sections of certain line bundles over $S$. All of this makes Kronheimer's hyperk\"ahler structure of an orbit somewhat more explicit: there are no differential equations left to solve. In particular, following up on ideas of Hitchin \cite{Hit2}, we give a formula for a K\"ahler potential of Kronheimer's hyperk\"ahler metric in terms of the theta function of $S$. We also interpret the hypercomplex structure of an open dense subset of the orbit as a natural structure on a $PU(k)$ bundle over a subset of the Jacobian of $S$ and identify the twistor lines as intersections of quadrics. We remark that a different characterisation of twistor lines has been given by Santa-Cruz \cite{SC}.

\tableofcontents

\section{Line bundles and flows on spectral curves\label{lb}}

\subsection{Line bundles and matricial polynomials}

In what follows $T$ denotes the total space of the line bundle $\sO(2)$ on $\oP^1$ ($T\simeq T\oP^1$), $\pi:T\rightarrow \oP^1$ is the
projection, $\zeta$ is the affine coordinate on $\oP^1$ and $\eta$ is the fibre coordinate on $T$. In other words $T$ is obtained by gluing
two copies of $\oC^2$ with coordinates $(\zeta,\eta)$ and $(\tilde{\zeta},\tilde{\eta})$ via:
$$ \tilde{\zeta}=\zeta^{-1}, \quad \tilde{\eta}=\eta/\zeta^2.$$
We denote the corresponding two open subsets of $T$ by $U_0$ and $U_\infty$.

Let $S$ be an algebraic curve in the linear system $\sO(2k)$, i.e. over $\zeta\neq \infty$ $S$ is defined by the equation
\begin{equation} P(\zeta,\eta)= \eta^k+a_1(\zeta)\eta^{k-1}+\cdots +a_{k-1}(\zeta)\eta+ a_k(\zeta)=0,\label{S}\end{equation}
where  $a_i(\zeta)$ is a polynomial of degree $2i$. $S$ can be singular or non-reduced.
\par
We recall the following facts (see, e.g., \cite{Hit,AHH}):
\begin{proposition} The group $H^1(T,\sO_T)$ (i.e. line bundles on $T$ with zero first Chern class) is generated by  $\eta^i\zeta^{-j}$, $i>0$, $0<j<2i$. The corresponding line bundles have transition functions $\exp(\eta^i\zeta^{-j})$ from $U_0$ to $U_\infty$.\hfill $\Box$\label{T}\end{proposition}
\begin{proposition} The natural map $H^1(T,\sO_T)\rightarrow H^1(S,\sO_S)$ is a surjection, i.e. $H^1(S,\sO_S)$ is generated by $\eta^i\zeta^{-j}$, $0<i\leq k-1$, $0<j<2i$.\hfill $\Box$\label{all}\end{proposition}

Thus, the (arithmetic) genus of $S$ is $g=(k-1)^2$. For a smooth $S$, the last proposition describes line bundles of degree $0$ on $S$.
In general, by a line bundle we mean an invertible sheaf. Its degree is defined  as its Euler characteristic plus $g-1$. The theta divisor $\Theta$ is  the set of line bundles  of degree $g-1$ which have a non-zero section.
\par
Let $\sO_T(i)$ denote the pull-back of $\sO(i)$ to $T$ via $\pi:T\rightarrow \oP^1$. If $E$ is a sheaf on $T$ we denote by $E(i)$ the sheaf
$E\otimes \sO_T(i)$ and similarly for sheaves on $S$. In particular, $\pi^\ast \sO$ is identified with $\sO_S$.
\par
If $F$ is a line bundle of degree $0$ on $S$, determined by a cocycle $q\in H^1(T,\sO_T)$, and $s\in H^0\bigl(S, F(i)\bigr)$, then we denote by
$s_0,s_\infty$ the representation of $s$ in the trivialisation $U_0,U_\infty$, i.e.:
\begin{equation} s_\infty(\zeta,\eta)=\frac{e^q}{\zeta^i}s_0(\zeta, \eta).\label{represent}\end{equation}

 We recall the following theorem of Beauville \cite{Beau}:
%%%%
\begin{theorem} There is a $1-1$ correspondence between the affine Jacobian $J^{g-1}-\Theta$ of line bundles of degree $g-1$ on $S$ and $GL(k,\cx)$-conjugacy classes of $\gl(k,\cx)$-valued polynomials $A(\zeta)=A_0+A_1\zeta+A_2\zeta^2$ such that $A(\zeta)$ is regular for every $\zeta$ and the characteristic polynomial of $A(\zeta)$ is \eqref{S}.\hfill ${\Box}$\label{Beauville} \end{theorem}

The correspondence is given by associating to a line bundle $E$ on $S$ its direct image $V=\pi_\ast E$, which has a structure of a $\pi_\ast \sO$-module. This is the same as a homomorphism $A:V\rightarrow V(2)$ which satisfies \eqref{S}. The condition $E\in J^{g-1}-\Theta$ is equivalent to $H^0(S,E)=H^1(S,E)=0$ and, hence, to $H^0(\oP^1,V)=H^1(\oP^1,V)=0$, i.e. $V=\bigoplus \sO(-1)$. Thus, we can interpret $A$ as a matricial polynomial precisely when $E\in J^{g-1}-\Theta$.

Somewhat more explicitly, the correspondence is seen from the exact sequence
\begin{equation} 0\rightarrow \sO_T(-2)^{\oplus k}\rightarrow \sO_T^{\oplus k}\rightarrow E(1)\rightarrow 0, \label{bundle}\end{equation}
where the first map is given by $\eta\cdot 1-A(\zeta)$ and $E(1)$ is viewed as a sheaf on $T$ supported on $S$. The inverse map is defined by the commuting diagram
\begin{equation}\begin{CD} H^0\bigl(S,E(1)\bigr) @>>> H^0\bigl(D_{\zeta}, E(1)\bigr)\\ @V \tilde{A}(\zeta) VV  @VV \cdot \eta V \\ H^0\bigl(S,E(1)\bigr) @>>> H^0\bigl(D_{\zeta}, E(1)\bigr), \end{CD}
\label{endom}\end{equation} where $D_{\zeta}$ is the divisor consisting of points of $S$ which lie above $\zeta$ (counting multiplicities).
That the endomorphism $\tilde{A}(\zeta)$ is quadratic in $\zeta$ is proved e.g. in \cite{AHH}. Observe that if $D_{\zeta_0}$ consists of
$k$ distinct points $p_1,\dots,p_k$ and if $\psi^1,\dots \psi^k$ is a basis of $H^0\bigl(S,E(1)\bigr)$, then $\tilde{A}(\zeta_0)$ in this basis is
\begin{equation} A(\zeta_0)=\left[ \psi^j(p_i)\right]^{-1} \diag\bigl(\eta(p_1),\dots,\eta(p_k)\bigr)\left[ \psi^j(p_i)\right],\label{conjugate} \end{equation}
where $\left[ \psi^j(p_i)\right]$ is a matrix with rows labelled by $i$ and columns by $j$.
%%%% extend to compactified Jacobian
%%%%
\begin{remark} For a singular curve $S$, Beauville's correspondence most likely extends to $\overline{J^{g-1}}-\overline{\Theta}$, where $\overline{J^{g-1}}$ is the compactified Jacobian  in the sense of \cite{Alex}. It seems to us that this is essentially proved in \cite{AHH}.\label{comp}\end{remark}

Let $K$ be the canonical (or dualising) sheaf of $S$. We have $K\simeq  \sO_S(2k-4)$. If $E$ belongs to $J^{g-1} -\Theta$, then so does $E^\ast \otimes K$ and:
\begin{proposition} Let $A(\zeta)$ be the quadratic matricial polynomial corresponding to $E\in J^{g-1}-\Theta$. Then $A(\zeta)^T$ corresponds to $E^\ast \otimes K$. \label{canonical}\end{proposition}
\begin{proof} For a finite morphism $f:X\rightarrow Y$ of projective schemes and coherent $\sO_X$-module $E$ and $\sO_Y$-module $G$ we have a natural isomorphism \cite[Ex. III.6.10]{Hart}:
$$ f_\ast \text{\it Hom}_X\bigl(E,f^!G)\stackrel{\sim}{\rightarrow}\text{\it Hom}_Y(f_\ast E,G),\enskip \text{where}\enskip f^!G=\text{\it Hom}_Y(f_\ast \sO_Y,G).$$
On the other hand, denoting by $K_X$ and $K_Y$ the dualising sheaves of $X$ and $Y$, we have from \cite[Ex. III.7.2]{Hart} that $f^!K_Y=K_X$. Therefore, in our situation, where $X=S$, $Y=\oP^1$, $f=\pi$, $\pi_\ast E=V$, we have $\pi_\ast \bigl(E^\ast \otimes K\bigr)\simeq V^\ast(-2)$.
\end{proof}
%%%%%%%%%%%%%
In particular, theta-characteristics outside $\Theta$ correspond to symmetric matricial polynomials.

\subsection{Real structure\label{realsect}} The space $T$ is equipped with a real structure (i.e. an antiholomorphic involution) $\tau$ defined by
\begin{equation} \zeta\mapsto -\frac{1}{\bar{\zeta}},\quad \eta\mapsto -\frac{\bar{\eta}}{\bar{\zeta}^2}. \label{antipodal}\end{equation}
%%%
Suppose that $S$ is real, i.e. invariant under $\tau$. Then $\tau$ induces an antiholomorphic involution $\sigma$ on $\Pic S$ as follows. Let $E$ be a line bundle on $S$ trivialised in a cover $\{U_\alpha\}_{\alpha\in A}$ with transition functions $g_{\alpha\beta}(\zeta,\eta)$ from $U_\alpha$ to $U_\beta$. Then $\sigma(E)$ is trivialised in the cover $\bigl\{\tau(U_\alpha)\bigr\}_{\alpha\in A}$ with transition functions
$$ \overline{g_{\alpha\beta}\bigl(\tau(\zeta,\eta)\bigr)}, $$
from $\tau(U_\alpha)$ to $\tau(U_\beta)$.  Observe that $\sigma(E)= \overline{\tau^\ast E}$
where ``bar" means taking the opposite complex structure. This map does not change the degree of $E$ and preserves line bundles $ \sO_S(i)$.
As there is a corresponding map on sections
\begin{equation} \sigma: s\mapsto \overline{\tau^\ast s},\label{sigma0}\end{equation}
it is clear that $J^{g-1} -\Theta$ is invariant under this map. Its effect on matricial polynomials is given by:
\begin{lemma} Let $A(\zeta)=A_0+A_1\zeta+A_2\zeta^2$ be the quadratic matricial polynomial corresponding to $E\in J^{g-1}-\Theta$. Then $-\overline{A_2}+\overline{A_1}\zeta-\overline{A_0}\zeta^2$ corresponds to $\sigma(E)$. \label{tauE}\end{lemma}
\begin{proof} From \eqref{endom}, using \eqref{sigma0}, the polynomial  $A^\tau(\zeta)$ corresponding to $\sigma(E)$ satisfies:  $$A^\tau(-1/\bar{\zeta})=-\overline{A(\zeta)}/\bar{\zeta}^2,$$
from which the statement follows.\end{proof}
%%%%%%%%%%%%%%%%%%%%%%%%%%%%
As we are interested in the hermitian conjugation of $A(\zeta)$, and not in the complex one, Proposition \ref{canonical} leads us to adopt the following definition:
\begin{definition} A line bundle $E$ of degree $d=ik$ on $S$, $i\in \oZ$, is {\em real} if $E\simeq \sigma(E)^{\ast}\otimes \sO_S(2i)$. We denote the corresponding
subspace of $J^{d}(S)$ by $J_{\oR}^{d}(S)$ and, in the case $d=g-1=(k-2)k$ write $\Theta_\oR$ for $\Theta\cap J_{\oR}^{g-1}$.\label{real-def}\end{definition}
%%%%%%%%%%%%%
Proposition \ref{T} and \ref{all} imply that real line bundles of degree $0$ have
transition functions $\exp q(\zeta,\eta)$, where $q$ satisfies:
$$\overline{q(\tau(\zeta,\eta))}=q(\zeta,\eta).$$
In general, we have:
\begin{lemma} A line bundle $E$ of degree $d=ik$, $i\in \oZ$, on $S$ is real if and only if it is of the form $E=F(i)$, where $F$ is a real line bundle of degree $0$.\hfill{$\Box$}\label{real-i}\end{lemma}
\begin{remark} This lemma is empty for a smooth $S$. However, on a reducible curve $S=\bigcup_{j=1}^p S_j$, $S_j\in |\sO(2k_j|$, it implies that a real bundle of degree $ik$ has the multi-degree $(ik_1,\dots, ik_p)$. Thus, the reality condition  picks out one component of $J^d (S)$.\label{multi}\end{remark}

For bundles of degree $g-1$ we conclude:
%%%
\begin{proposition} There is a $1-1$ correspondence between $J_{\oR}^{g-1}-\Theta_\oR$ and conjugacy classes of matrix-valued polynomials $A(\zeta)$ as in Theorem \ref{Beauville} such that there exists a hermitian $h\in GL(k,\cx)$ with
\begin{equation} hA_0h^{-1}=-A_2^\ast,\quad hA_1h^{-1}=A_1^\ast,\quad hA_2h^{-1}=-A_0^\ast.\label{h}\end{equation}
\label{realbundles}\end{proposition}
\begin{proof} From Proposition \ref{canonical} and Lemma \ref{tauE}, we know that a bundle in $J^{g-1}-\Theta$ is real if and only if $A(\zeta)=A_0+A_1\zeta+A_2\zeta^2$ is conjugate to $-A_2^\ast+A_1^\ast\zeta-A_0^\ast\zeta^2$. Let $h$ be the matrix realising this conjugation, i.e. $h$ satisfies \eqref{h}. It follows that $p=(h^\ast)^{-1}h$ centralises $A(\zeta)$. If $p$ does not belong to the centre of $GL(k,\cx)$, then $A(\zeta)$ belongs to the reductive subalgebra $Z(p)$ and $\Ker(\eta-A(\zeta))$ is not everywhere one-dimensional. This contradicts the assumption that $A(\zeta)$ corresponds to a line bundle. Thus $h$ is hermitian up to an irrelevant central factor.
\end{proof}

\begin{example} This example illustrates Remark \ref{multi} in view of the above proposition. Consider a curve $S$ consisting of two rational components intersecting in a pair of points (i.e. the polynomial \eqref{S} is the product of two polynomials linear in $\eta$). Then $k=2$ and $g=1$. The affine Jacobian  $J^{g-1}-\Theta=J^0-\Theta$ has $3$ components: the line bundles of multi-degree $(0,0)$, apart from $ \sO_S$, and all line bundles of multi-degrees $(1,-1)$ and $(-1,1)$. The latter two components correspond to matricial polynomials which are  upper-triangular for all $\zeta$ (up to conjugation by $GL(2,\cx)$) with the diagonal part distinguishing $(1,-1)$ from $(-1,1)$. The last proposition shows again that there are no real bundles in these components.
\end{example}

\subsection{Hermitian metrics\label{hermitian}}

\begin{definition} A line bundle of degree $g-1$ on $S$ is called {\em definite} if it is in $J_{\oR}^{g-1}-\Theta_\oR$ and the matrix $h$ in \eqref{h}
can be chosen to be positive-definite. The space of definite line bundles on $S$ will be denoted by
$J^{g-1}_+$.\label{pos}\end{definition}

We easily conclude that there is a 1-1 correspondence between $J_+^{g-1}$ and $U(k)$-conjugacy classes of matrix-valued polynomials
$A(\zeta)$ as in Theorem \ref{Beauville} which in addition satisfy
\begin{equation} A_2=-A_0^\ast,\quad A_1=A_1^\ast. \label{+} \end{equation}

Definite line bundles have also the following interpretation (cf. \cite{Hit}):
\par
For $E=F(k-2)\in J_\oR^{g-1}$ the real structure induces an antiholomorphic isomorphism
\begin{equation}\sigma: H^0\bigl(S, F(k-1)\bigr)\longrightarrow H^0\bigl(S, F^{\ast}(k-1)\bigr),
\label{sigma}\end{equation} via the map \eqref{sigma0}. Thus, for   $v,w\in H^0\bigl(S,
F(k-1)\bigr)$,  $v\sigma(w)$ is a section of $\sO_S(2k-2)$ and so it can be uniquely written \cite{Hit,AHH} as
\begin{equation} c_0\eta^{k-1}+ c_1(\zeta)\eta^{k-2}+\dots +c_{k}(\zeta), \label{c_0}
\end{equation}
where the degree of $c_i$ is $2i$. Following Hitchin \cite{Hit}, we define a hermitian form on $H^0\bigl(S, F(k-1)\bigr)$ by
\begin{equation} \langle v,w\rangle=c_0.\label{form}\end{equation}
%%%%%%%%%%%%%%%%%
The following fact can be deduced from \cite{Hit}:
\begin{proposition} A line bundle $E=F(k-2)\in J_{\oR}^{g-1}-\Theta_\oR$ is in $J_+^{g-1}$ if and only if the above form on $H^0\bigl(S, F(k-1)\bigr)$ is definite.\hfill $\Box$
\end{proposition}

Because of this, we shall also refer to $F(k-1)$ as definite.

\subsection{Flows\label{flows}} If we fix a tangent direction on $J^{g-1}(S)$, i.e. an element $q$ of $H^1(S,\sO_S)$, then the linear flow of line bundles on
$J^{g-1}(S)$ corresponds to a flow of matricial polynomials (modulo the action of $GL(n,\cx)$). We shall be interested only in the flow
corresponding to $[\eta/\zeta] \in H^1(S,\sO_S)$. Following the tradition, we denote by $L^t$ the line bundle on $T$ with transition
function $\exp(-t\eta/\zeta)$ from $U_0$ to $U_\infty$.
\par
For any line bundle $F$ of degree $0$ on $S$ we denote by $F_t$ the line bundle $F\otimes L^t$. We consider the flow $F_t(k-2)$ on
$J^{g-1}(S)$. Even if $F=F_0$ is in the theta divisor, this flow transports one immediately outside $\Theta$, and so we obtain a flow of
endomorphisms of $V_t=H^0\bigl(S,F_t(k-1)\bigr)$. These vector spaces have dimension $k$ as long as $F_t(k-2)\not\in \Theta$. We obtain an
endomorphism $\tilde{A}(\zeta)$ of $V_t$ as equal to multiplication by $\eta$ on $H^0(S\cap \pi^{-1}(\zeta), F_t(k-1)\bigr)$, where
$\pi:T\rightarrow \oP^1$ is the projection.
\par
To obtain a flow of matricial polynomials one has to trivialise the vector bundle $V$ over $\oR$ (whose fibre at $t$ is $V_t$). This is a
matter of choosing a connection. If we choose the connection $\nabla^0$ defined by evaluating sections at $\zeta=0$ (in the trivialisation
$U_0,U_\infty$), then the corresponding matricial polynomial $A(t,\zeta)=A_0(t)+A_1(t)\zeta+A_2(t)\zeta^2$ satisfies \cite{Hit,AHH}
$$ \frac{d}{dt}A(t,\zeta)=\left[A(t,\zeta), A_2(t)\zeta\right].$$
%%%%%%%%%%%%%%%%%%%%%%
As mentioned above, if $F$ is a real bundle, then $V$ has a natural hermitian metric \eqref{form} (possibly indefinite). The above
connection is not metric, i.e. it does not preserve the form  \eqref{form}. Hitchin \cite{Hit} has shown that the connection
$\nabla=\nabla^0+\frac{1}{2}A_1(t)dt$ is metric and that, in a $\nabla$-parallel basis, the resulting $A(t,\zeta)$ satisfies
$$\frac{d}{dt}A(t,\zeta)=\left[A(t,\zeta), A_1(t)/2+A_2(t)\zeta\right].$$
If the bundle $F(k-1)$ is positive-definite, then so are all $F_t(k-1)$. If the basis of sections is, in addition, unitary, then the
polynomials $A(t,\zeta)$ satisfy the reality condition \eqref{+}. If we write $A_0(t)=T_2(t)+iT_3(t)$ and $A_1(t)=2iT_1(t)$ for
skew-hermitian $T_i(t)$, then these matrices satisfy the Nahm equations:
\begin{equation}\dot{T}_i+\frac{1}{2}\sum_{j,k=1,2,3}\epsilon_{ijk}[T_j,T_k]=0\;,\;\;\;\;
i=1,2,3.\label{Nahm0}\end{equation}

\section{A unitary basis\label{basis}}

For a definite line bundle $F(k-1)$ on a real curve $S$, we can find an explicit basis of sections, unitary with respect to the form
\eqref{form}. We begin by reformulating the equation \eqref{form}.

Let $s,s^\prime$ be two sections of $F(k-1)$ on $S$. The form $\langle s, s^\prime\rangle$ is given by computing the section
$Z=s\sigma(s^\prime)$ of $\sO(2k-2)$ on $S$. Writing
$$ Z(\zeta,\eta)=c_0\eta^{k-1}+ c_1(\zeta)\eta^{k-2}+\dots +c_{k}(\zeta)$$
on $S$, we have $\langle s, s^\prime\rangle=c_0$. If $P(\zeta,\eta)=0$ is the equation defining $S$, then for any $\zeta_0$, such that
$S\cap \pi^{-1}(\zeta_0)$ consists of distinct points, we have
$$ c_0= \sum_{(\zeta_0,\eta)\in S} \Res\frac{Z(\zeta_0,\eta)}{P(\zeta_0,\eta)}.$$
Thus, if we write $(\zeta_0,\eta_1),\dots, (\zeta_0,\eta_k)$ for the points of $S$ lying over $\zeta_0$, then we have
%%%%%%%%%%%%%%%%%%
\begin{equation} \langle s, s^\prime\rangle=\sum_{i=1}^k\frac{s(\zeta_0,\eta_i)\cdot\sigma(s^\prime)(\zeta_0,\eta_i)}{\prod_{j\neq i}
\bigl(\eta_i-\eta_j\bigr)}.\label{c0}\end{equation}
%%%%%%%%%%%%%%%%%
Therefore, one can compute $\langle s, s^\prime\rangle$ from the values of the sections at two fibres of $S$ over two antipodal points of
$\oP^1$ (as long as the fibres do not have multiple points).
\par
Without loss of generality, we shall take $0$ and $\infty$ as a pair of antipodal points in $\oP^1$. We assume that the fibre over both
points consists of distinct points, which denote by $0_i$ and $\infty_i$, $i=1,\dots k$, so that $\tau(0_i)=\infty_i$. We have
%%%%%%%%%%%%%%%%%%%%%%%%%
\begin{lemma} Let $F(k-1)$ be a definite line bundle on $S$. Let $K\subset \{1,\dots,k\}$ and let $D=\{0_i;i\in K\}\cup \{\infty_j;j\not\in K\}$.
Then $F(k-1)[-D]\not\in \Theta$. \label{D}\end{lemma}
\begin{proof} Let $s$ be a section of $F(k-1)[-D]$. Then, the formula \eqref{c0}, with $\zeta_0=0$, implies that $\|s\|^2=0$, so $s\equiv
0$.\end{proof}
%%%%%%%%%%%%%%
Hence, for any $i=1,\dots,k$, there is a one-dimensional space $V_i$ of sections vanishing at each $\infty_j$ for $j<i$ and each $0_j$ for
$j>i$. A section in $V_i$ can be found explicitly in terms of the theta function of $S$ (as in \cite{AHH}), which we shall do later for a
reducible curve. These sections are automatically mutually orthogonal:
%%%%%%%%%%%%%%%%%
\begin{proposition} The subspaces $V_i$ of $H^0\bigl(S,F(k-1)\bigr)$ are mutually orthogonal in the metric \eqref{form}. \label{V_i}
\end{proposition}
\begin{proof} Let $s\in V_i$ and $s^\prime\in V_j$, $i\neq j$. Then $s\sigma(s^\prime)$ vanishes at all $\infty_n$ with $n<i$ or $n>j$, and at
all $0_n$ with $n>i$ or $n<j$. Thus $s\sigma(s^\prime)$ vanishes at all $\infty_n$ or at all $0_n$. The formula \eqref{c0} shows that
$\langle s, s^\prime\rangle=0$.
\end{proof}

\begin{remark} The formula \eqref{conjugate} shows that the matricial polynomial $A(\zeta)$, corresponding to multiplication by $\eta$ in
a basis of sections $s^i\in V_i$, is upper-triangular at $\zeta=0$ and lower-triangular at $\zeta=\infty$.\label{vartheta}
\end{remark}

\section{Reducible algebraic curves}

In this section, $S$ does not have to be, a priori, in $|\sO(2k)|$. We consider a curve $S$ whose consisting of  $k$ copies $S_1,\dots,S_k$
of $\oP^1$, such that each pair of the $S_i$ intersects in two points and such that all intersection points are distinct. The genus $g$ of
$S$ is $(k-1)^2$. We assume, in addition, that $S$ comes with a map $\pi:S\rightarrow \oP^1$ whose restriction to each $S_i$ is an
isomorphism.
\par
For any $i\neq j$, we denote by $a_{ij},a_{ji}$ the image under $\pi$ of the two intersection points of $S_i$ and $S_j$. Since all
singularities of $S$ are ordinary double points, we have: $a_{ij}\neq a_{ik}$ and $a_{ji}\neq a_{ki}$ if $j\neq k$. We shall also assume
that $a_{ij}\neq \infty $, $i, j=1,\dots,k$, $i\neq j$.
\par
Let $p:\tilde{S}\rightarrow S$ be the normalisation. We denote by $\phi_i:\oP^1\rightarrow S_i$ the inverse of $\pi\circ p|_{S_i}$.

\subsection{Line bundles\label{line}} We denote by $\Pic^0 S$ the group of line bundles of multi-degree $0=(0,\dots,0)$. Such a line bundle $F$  on $S$ gives $k$ copies $\sO_i\simeq \sO$ of the trivial bundle on each $S_i$, $i=1,\dots,k$, and hence $F$
is identified by the matching conditions at each $a_{ij}$, i.e. by an invertible automorphism $\lambda_{ij}\in \oC^\ast$ between the fibre
of $\sO_i$ at $\phi_i(a_{ij})$ and the one of $\sO_j$ at $\phi_j(a_{ij})$. A change of the isomorphisms $\sO_i\simeq \sO$ does not change
the bundle $F$, and hence, we identify $\Pic^0 S$ with $\bigl(\oC^\ast\bigr)^{k(k-1)}/ \bigl(\oC^\ast\bigr)^{k}$ where the action of
$\bigl(\oC^\ast\bigr)^{k}$ is given by
\begin{equation} (z_1,\dots,z_k)\cdot (\lambda_{ij})=\bigl( z_iz_j^{-1}\lambda_{ij}\bigr).\label{action}\end{equation}
Observe that the diagonal $\oC^\ast$ acts trivially and, hence, $\Pic^0 S\simeq\bigl(\oC^\ast\bigr)^{(k-1)^2}$. We denote by
$[\lambda_{ij}]$ the point in $\Pic^0 S$ corresponding to a $(\lambda_{ij})\in\bigl(\oC^\ast\bigr)^{k(k-1)}$.

\subsection{The Abel map} We define the Jacobian of $S$ (the range of abelian integrals) as $\Jac S=\bigl(\cx^\ast\bigr)^g\simeq \Pic^0 S$.
We consider line bundles $F(n)=F\otimes \pi^\ast \sO(n)$ on $S$, where $F\in \Pic^0 S$ and $n\geq 0$, i.e. line bundles of
multi-degree $(n,\dots,n)$. We denote by $S^{(n)}$ the set  of effective divisors of multi-degree $(n,\dots,n)$, i.e.
$$ S^{(n)}=\left\{(x_{i}^m)_{\stackrel{\scriptstyle i=1,\dots,k}{\scriptstyle m=1,\dots,n}} ;  \enskip x_l^m\in S_l, x_l^m\neq a_{lj}, a_{jl}\right\}$$
(we omit the reference to the isomorphisms $\phi_i$; strictly speaking we have $\pi\circ p(x_l^m)\neq a_{lj}, a_{jl}$).
%%%%%%%%%%%%%%%%%%%%%555555
We define polynomials $P_i(\zeta)=\prod_{m=1}^n(\zeta-x_i^m)$ (if any $x_i^m$ is $\infty$, then that factor is omitted). These give a
section of $F(n)$, where $F$ is determined by $[\lambda_{ij}]$, if and only if $P_j(a_{ij})=\lambda_{ij}P_i(a_{ij})$ for every $i,j$, $i\neq
j$. Therefore, we can define the Abel map $A_n=A_n(\infty,\cdot):S^{(n)}\rightarrow \Jac S$ by
$$ A_n\bigl((x_{i}^m)\bigr)=[\lambda_{ij}]\enskip\text{where}\enskip \lambda_{ij}=\frac{\prod_{m=1}^n(a_{ij}-x_j^m)}{\prod_{m=1}^n(a_{ij}-x_i^m)}.$$

\subsection{Theta function and theta divisor} In $J^{g-1} S$ we again consider only the component of line bundles of multi-degree $(k-2,\dots,k-2)$, i.e. line bundles of the form  $F(k-2)$, where $F\in \Pic^0 S$. In this component, the theta divisor is the image of
$A_{k-2}$, i.e. the set of line bundles $F(k-2)$  which have a section. Concretely, $\Theta$ is the set of equivalence classes of
$[\lambda_{ij}]$ of $(\lambda_{ij})\in \bigl(\oC^\ast\bigr)^{k(k-1)}$ for which there exist polynomials $Q_1(\zeta),\dots, Q_k(\zeta)$ of
degree $\leq k-2$, which solve the equations
\begin{equation} Q_j(a_{ij})=\lambda_{ij}Q_i(a_{ij}),\quad i,j=1,\dots,k,\enskip i\neq j. \label{match1}\end{equation}
%%%%%%%%%%%%%%%%%%%%%%%%%
These equations can be written as $k(k-1)$ linear equations $\Xi Q=0$, where $\Xi$ is a matrix depending on the $a_{ij}$ and on the
$\lambda_{ij}$ and $Q$ is the vector of the coefficients of polynomials $Q_i$. We shall see shortly that the determinant of $\Xi$ is
invariant under \eqref{action}, and hence, the function
\begin{equation} \vartheta(\lambda_{ij})=\det\Xi\label{theta0}\end{equation}
descends to a function on $\Jac S$.  This function is the analogue of the theta function (since there are no $B$-cycles on $S$, the theta
function lives on $\Jac S$ rather than just on its universal cover). The set $\Theta$ of zeros of $\vartheta$ is the theta divisor of line
bundles $F$ such that $F(k-2)$ admits a non-zero section. We use the notation $\vartheta$, instead of $\theta$, to emphasise that our theta
function is already ``translated by the Riemann constant".

For many purposes, it is better to define the theta function on a bigger space, given by replacing \eqref{match1} with homogeneous
equations
\begin{equation}\lambda_{ij}Q_i(a_{ij})+\mu_{ij}Q_j(a_{ij})=0.\label{hom} \end{equation}
%%%%%%%%%%%%%%%%%%%%%%%%%%%%%%%%%%%
Thus we consider the space $\bigl(\oC^\ast \times \oC^\ast\bigr)^{k(k-1)}$ of pairs $(\lambda_{ij},\mu_{ij})$, which is a fibration over
$\Jac S$ under the map which sends $(\lambda_{ij},\mu_{ij})$ to the orbit (under \eqref{action}) of the quotient $-\lambda_{ij}/\mu_{ij}$.
To define a determinant we have to make \eqref{hom} into a matrix $\Xi$. The rows of $\Xi$ are the equations \eqref{hom} for every $(i,j)$ with $(i,j)$ ordered lexicographically. The columns of $\Xi$ are numbered by pairs $(m,n)$, $m=1,\dots,k$, $n=0,\dots,k-2$, with each $(m,n)$ corresponding  to the $n$-th
coefficient of $Q_m(\zeta)=b_0+\dots+b_{k-2}\zeta^{k-2}$. The ordering on $(m,n)$ is again lexicographical.

With these preliminaries, we define:
\begin{definition}   The (extended) theta function of the curve $S$ is the map
$\vartheta:\bigl(\oC^\ast \times \oC^\ast\bigr)^{k(k-1)}\rightarrow \oC,$ defined by
$$ \vartheta(\lambda_{ij},\mu_{ij})=\det \Xi,$$
where $\Xi$ is the matrix given by the left-hand side of
\eqref{hom}. \end{definition}

An immediate advantage of defining $\vartheta$ this way is that it depends on a particular labelling of intersection points only up to a sign.

The projection of the set of zeros of $\vartheta$ to $\Jac S$ is the theta divisor.
 It is clear that $\vartheta$ is a homogeneous polynomial in $\lambda_{ij},\mu_{ij}$, but not all degrees occur. We adopt the following
definition:
%%%%%%%%%%%%%%%%%%%%
\begin{definition} Let $P=\{(i,j); i,j=1,\dots,k,\enskip i\neq j\}$. A subset $L$ of $P$ is called {\em regular}, if, for every $i$, there are as many elements of $L$ with first coordinate $i$ as with second coordinate equal to $i$. In other words, for every $i=1,\dots, k$, the set
\begin{equation} L_i=\{(i,j)\in L\}\cup \{(m,i)\not\in L\}\label{L_i}
\end{equation}
has cardinality $k-1$. The set of regular subsets of $P$ will be denoted by $\sR(P)$.\end{definition}
%%%%%%%%%%%%%%%%
%%%%%%%%%%%%%%%%%%%%
\begin{theorem} The following formula for the theta function $\vartheta:\bigl(\oC^\ast \times \oC^\ast\bigr)^{k(k-1)}\rightarrow \oC$ of the curve $S$ holds true:
\begin{equation} \vartheta(\lambda_{ij},\mu_{ij})=\sum_{L\in \sR(P)} a_L\cdot \prod_{(i,j)\in L} \lambda_{ij}\prod_{(i,j)\not\in L} \mu_{ij},\label{theta}\end{equation}
where each coefficient $a_L$ is given by
$$ a_L=(-1)^s\prod_{m=1}^k V(L_m);$$
$V(L_m)$ is a Vandermonde determinant of $\{a_{ij};(i,j)\in L_m\}$ (for some ordering) and $s$ is the sign of the permutation of $P$ defined as $(L_1,\dots,L_k)$, where the ordering in each $L_m$ is the one used to define $V(L_m)$.
\label{Theta}\end{theorem}
%%%%
\begin{remark} Mumford \cite{mum} gives a similar formula in the case of a singular curve whose normalisation is $\oP^1$.\end{remark}
\begin{proof}
The matrix $\Xi$ has size $k(k-1)\times k(k-1)$. Its rows are numbered by pairs $(i,j)$, $i\neq j$, corresponding to the intersection
points $a_{ij}$. Its columns are numbered by pairs $(m,n)$, $m=1,\dots,k$, $n=0,\dots,k-2$, with each $(m,n)$ corresponding  to the $n$-th
coefficient of $Q_m(\zeta)=b_0+\dots+b_{k-2}\zeta^{k-2}$. Thus the entry of $\Xi$ with coordinates $(i,j),(m,n)$ is $\lambda_{ij}a_{ij}^n$
if $m=i$, $\mu_{ij}a_{ij}^n$ if $m=j$,  and it is zero otherwise. To compute the determinant of $\Xi$, we therefore have to sum up over
$1-1$ maps $\sigma$ from $P$ to $\{(m,n); m=1,\dots,k, n=0,\dots,k-2\}$ such that $\sigma(i,j)=(i,n)$ or $\sigma(i,j)=(j,n)$, for every
$(i,j)$. Observe that such a map corresponds to giving a subset $L$ of $P$ and a pair of maps $\rho:L\rightarrow \{0,\dots,k-2\}$,
${\rho}^\prime:L^\prime\rightarrow \{0,\dots,k-2\}$ ($L^\prime$ is the complement of $L$), so that $\sigma(i,j)$ is $(i,\rho(i,j))$ or
$(j,\rho^\prime (i,j))$. The bijectivity amounts to the following conditions: (1) $\rho(i,\cdot)$ and $\rho^{\prime}(\cdot, j)$ are
injective for every $i,j$; (2) $\rho(i,j)\neq \rho^{\prime}(m, i)$ if $(i,j)\in L$ and $(m, i)\in L^\prime$. These imply that, for every
$i$, we have an injective map $\sigma_i$ from $L_i$ given by \eqref{L_i}  to $\{0,\dots,k-2\}$ given by $\sigma_i(i,j)=\rho(i,j)$ and
$\sigma_i(m,i)=\rho^{\prime}(m, i)$. The sets $L_i$ are disjoint and as every $(i,j)$ belongs to $L_i$ or $L_j$, their union is $P$. Since
each $\sigma_{i}$ is injective, each $L_i$ has cardinality at most $k-1$ and therefore exactly $k-1$. Thus $L$ is regular and each
$\sigma_i$ is a bijection. If we write $\Sigma(L)$ for the set of $k$-tuples $(\sigma_1,\dots,\sigma_k)$ such that each $\sigma_i$ is
a bijective mapping from $L_i$ onto $\{0,1,\dots,k-2\}$, then the determinant can be written as
$$\sum_{L\in \sR(P)}\sum_{\sigma\in \Sigma(L)} (-1)^{\sgn\sigma}
 \prod_{(i,j)\in L} \lambda_{ij}a_{ij}^{\sigma_i(j)}\prod_{(i,j)\not\in L} \mu_{ij}a_{ij}^{\sigma_j(i)},$$
 where $\sgn\sigma$ denotes the sign of the corresponding permutation of $P$. This proves the formula \eqref{theta}, and each coefficient can be evaluated by grouping together factors with coordinates in each $L_m$.
\end{proof}
%%%%%%%%%%%%55

We can recover the function \eqref{theta0} and similar functions on $\Jac S$ from the following observation
\begin{proposition} Let $p,q$ be two integers and define  a map
$$f:\bigl(\oC^\ast\bigr)^{k(k-1)}\rightarrow \bigl(\oC^\ast \times \oC^\ast\bigr)^{k(k-1)},\quad
 f(\lambda_{ij})=(\lambda_{ij}^p,\lambda_{ij}^q).$$ Then the map $\vartheta\circ f$ is invariant under the action \eqref{action} and induces a
 map $\vartheta_{p,q}:\Jac S\rightarrow \oC$.\label{p,q}\end{proposition}
 \begin{proof} This is a consequence of regularity of subsets $L$ (and their complements) occurring in \eqref{theta}: in any $L$, the number of
 $\lambda_{ij}$  with the first coordinate $m$ is the same as the number of those with second coordinate $m$. \end{proof}
%%%%%%%%%%%%%%%%%%
For nonnegative $p,q$, the function $\vartheta\circ f$, given in the statement, clearly extends to $\oC^{k(k-1)}$ and, hence, $\vartheta_{p,q}$ extends to the
algebro-geometric quotient of $\oC^{k(k-1)}$ by \eqref{action}. We continue to write $[\lambda_{ij}]$ for the orbit of $(\lambda_{ij})\in
\oC^{k(k-1)}$. From Theorem \ref{Theta}, we immediately obtain:
%%%%%%%%%%%%%%%%%%%
\begin{corollary} $\vartheta_{1,0}$ and $\vartheta_{0,1}$ do not vanish at the origin $[0]$, i.e. at the orbit of $0\in \oC^{k(k-1)}$.\hfill $\Box$\label{0}
\end{corollary}

\begin{remark} Geometrically, this corollary can be interpreted in terms of the compactified Jacobian of $S$ \cite{Beau1, Alex}. The compactified Jacobian
$\overline{J^{g-1}(S)}$ is stratified by Jacobians of partial normalisations of $S$ and the boundary of the theta divisor of $S$ is given by
the theta divisors of these normalisations. The point $[0]$ is the smallest stratum in $\overline{J^{g-1}(S)}$, i.e. it corresponds to the full normalisation
$\tilde{S}$. Of course, there are no line bundles of degree $g(\tilde{S})-1$ admitting nonzero sections.\label{compactified}\end{remark}

\subsection{The inverse mapping}
As is well known, the theta function may be used to invert the Abel mapping. We shall now do this explicitly for the curve $S$.
\par
We consider a point in $\Jac S$ given by $(\lambda_{ij},\mu_{ij})$. It determines a line bundle $F$ of multi-degree $0$. We consider the line
bundle $F(k-1)$ and its sections. A section of $F(k-1)$ consists of polynomials $Q_1(\zeta),\dots,Q_k(\zeta)$ which satisfy \eqref{hom}. We
wish to construct such sections explicitly.

We define functions
\begin{equation} s_{ij}(u)= \lambda_{ij}(a_{ij}-u), \enskip t_{ij}(u)= \mu_{ij}(a_{ij}-u),\label{u}\end{equation}
with the convention $s_{ij}(\infty)=\lambda_{ij}, t_{ij}(\infty)=\mu_{ij}$, and a function of $k$ variables
\begin{equation} \Lambda_F(u_1,\dots,u_k)= \vartheta\bigl(s_{ij}(u_i),t_{ij}(u_j)\bigr).\label{Lambda}\end{equation}
%%%%%%%%%%%%%%%%%%%%%%%%%
Observe that, since every $L$ in \eqref{theta} is regular, $\Lambda$ is a polynomial of degree $k-1$ in each variable.
%%%%%%%%%%%%%%%%%%%%%%%%%%%%%%%%
\begin{proposition} Let $p_i\in S_i$, $i=1,\dots,k$, be different from any intersection points. Let $F$ be a line bundle of multi-degree $0$ on
$S$. Then $F(k-1)[-p_1-\dots-p_k]\in \Theta$ if and only if $\Lambda_F\bigl(\pi(p_1),\dots,\pi(p_k)\bigr)=0$.\label{zeros}\end{proposition}
\begin{proof} Put $D=p_1+\dots +p_k$ and $y_i=\pi(p_i)$, $i=1,\dots,k$. Let $\bigl(Q_1(\zeta),\dots,Q_k(\zeta)\bigr)$  a section of $F(k-1)[-D]$, i.e.
each $Q_i$ has degree $k-1$,
$Q_i(y_i)=0$, and $\lambda_{ij}Q_i(a_{ij})+\mu_{ij}Q_j(a_{ij})=0$ for every $i\neq j$. We then have
$$ s_{ij}(u_i)\frac{Q_i(a_{ij})}{a_{ij}-u_i}+ t_{ij}(u_j)\frac{Q_j(a_{ij})}{a_{ij}-u_j}=0$$
as functions of $u_1,\dots,u_k$, for all $i\neq j$. The function $Q_{i}(\zeta)/(\zeta- u_i)$ becomes at $u_i=y_i$ a polynomial of degree
$k-2$ and, hence, the section $(Q_1,\dots,Q_k)$ is non-zero if and only if $\vartheta\bigl(s_{ij}(y_i),t_{ij}(y_j)\bigr)=0$.
\end{proof}

\begin{proposition} Let $F$ and $p_i$, $i=1,\dots,k$ be as in the previous proposition and suppose that $F(k-1)[-p_1-\dots-p_k]\not\in
\Theta$.  Let $s_l=\bigl(Q_1(\zeta),\dots,Q_k(\zeta)\bigr)$, $l=1,\dots,k$, be the (unique up to a constant multiple) section of $F(k-1)$
which vanishes at $p_j$ for all $j\neq l$. Then, up to a constant multiple, the polynomial $Q_l(\zeta)$ is given by
\begin{equation}Q_l(\zeta)=\Lambda_F\bigl(\pi(p_1),\dots,\underbrace{\zeta}_{l},\dots,\pi(p_k)\bigr). \label{Q_l}\end{equation}
\label{Ql}\end{proposition}
%%%%
\begin{proof} Let $x_1,\dots,x_{k-1}$ be the roots of $Q_l$ and $w_i=p(\phi_l(x_i))$, $i=1,\dots,k-1$, the corresponding points of $S_l$.
Thus the section $s_l$ vanishes at $p_1,\dots,\underbrace{w_i}_{l},\dots,p_k$ for every $i$ and, from the previous proposition, the
right-hand side of \eqref{Q_l} vanishes at $x_1,\dots,x_{k-1}$. By the remark before the statement of Proposition \ref{zeros}, the
right-hand side is a polynomial of degree $k-1$ and, so, it is either a constant multiple of $Q_l(\zeta)$ or it vanishes identically. The
assumptions imply, however, that $\Lambda_F\bigl(\pi(p_1),\dots,\pi(p_k)\bigr)\neq 0$ and the second possibility is thereby excluded.
\end{proof}

\begin{remark} To completely determine the section $s_l$  (i.e. to describe the $Q_j$, $j\neq l$), we now find one  root $x$ of $Q_l$ and apply
the formula \eqref{Q_l} with the index $l$ replaced by $j$, $j\neq l$, to points $y_1,\dots,\underbrace{x}_{l},\dots,y_k$ (where
$y_i=\pi(p_i)$).\end{remark}

\section{Reducible spectral curves\label{spec}}
We now relate the results of the previous section to those in sections \ref{lb} and \ref{basis}, i.e. we consider a fully reducible real
spectral curve $S$ given by the equation
\begin{equation}\prod_{i=1}^k \bigl(\eta-z_i-2x_i\zeta+\bar{z}_i\zeta^2\bigr)=0,
\label{reducible}\end{equation} where each $x_i$ is real. We assume that $S$ does not have multiple components, i.e. $(x_i,z_i)\neq
(x_j,z_j)$ for $i\neq j$.
\par
Let $S_i$, $i=1,\dots,k$, be the components, i.e.
\begin{equation} S_i=\{(\zeta,\eta); \eta=z_i+2x_i\zeta-\bar{z}_i\zeta^2\}, \quad i=1,\dots,k.
\label{S_i}\end{equation} Two curves $S_i$ and $S_j$ intersect in a pair of distinct points $p_{ij},p_{ji}=\tau(p_{ij})$,  where the
$\zeta$-coordinate $a_{ij}=\pi(p_{ij})$ of each $p_{ij}$ is given by
\begin{equation} a_{ij}=\frac{(x_i-x_j)+r_{ij}}{\bar{z}_i-\bar{z}_j}, \quad r_{ij}=\sqrt{(x_i-x_j)^2 +|z_i-z_j|^2}.\label{aij}\end{equation}
We assume that $S$ has only nodes, i.e. $p_{ij}\neq p_{mn}$ for $(i,j)\neq (m,n)$. By using the $SO(3)$-action on $T\oP^1$, we can
also assure that $a_{ij}\neq 0,\infty$ for any $i,j$.

\subsection{Line bundles}
 Let $F$ be a line bundle of multi-degree $0$ on $S$. From Propositions \ref{T} and \ref{all}, $F$ is trivialised on $U_0,U_\infty$ with the
transition function $\exp q(\zeta,\eta)$ where
\begin{equation} q(\zeta,\eta)=\sum_{i=1}^{k-1} \frac{\eta^i}{\zeta^i}q_i(\zeta),\quad q_i(\zeta)=\sum_{n=-i+1}^{i-1} d_{n,i}\zeta^n, \label{q}\end{equation}
for some complex numbers $d_{n,i}$. Moreover,  $F$ is real in the sense of Definition \ref{real-def} if $\overline{d_{n,i}}=(-1)^k d_{-n,i}$
for all $i,k$.
 For every $i=1,\dots,k$, we write $q\bigl(\zeta,z_i+2x_i\zeta-\bar{z}_i\zeta^2)$ as
$q^i_+(\zeta)+ q^i_-(\zeta)$, where $q^i_+$ (resp. $q^i_-$) contains all positive (resp. negative) powers of $\zeta$ and
$q^i_+(0)=q^i_-(0)$.
%%%%
The multiplication by $e^{-q^i_+(\zeta)}$ on $U_0$ and by $e^{q^i_-(\zeta)}$ on $U_\infty$ provides an isomorphism between the trivial
bundle on $S_i$ and $F|_{S_i}$. This gives us a canonical identification of $\Pic^0 S$ with $\bigl(\oC^\ast\bigr)^g$. In the notation of
the previous section, the bundle $F$ corresponds to the point $[\lambda_{ij}]\in \Jac S$, where
$\lambda_{ij}=e^{q_+^j(a_{ij})-q_+^i(a_{ij})}$, $i\neq j$. The property of $F$ being real is equivalent to
$[\lambda_{ij}]=[\overline{\lambda_{ji}}]$ for all $i\neq j$.

If we now consider the bundle $F(k-1)$ on $S$, then its section $s(\zeta)$ is given on each $S_i\cap U_0$ by
\begin{equation} s_0^i(\zeta)=e^{-q^i_+(\zeta)}Q_i(\zeta),
\label{comps}\end{equation} where the $Q_i$ are polynomials of degree $k-1$ satisfying the matching conditions. For the time-dependent line
bundle $F_t=F\otimes L^t$, where $L^t$ has the transition function $e^{-t\eta/\zeta}$, we have to replace $q(\eta,\zeta)$ with
$q(\zeta,\eta)-t\eta/\zeta$. Thus \eqref{comps} becomes:
\begin{equation}s_0^i(t,\zeta)=e^{t(x_i-\bar{z}_i\zeta)}e^{-q^i_+(\zeta)} Q_i(t,\zeta).\label{comps1}\end{equation}
Moreover, we observe that the line bundle $L^t$ corresponds to the point \begin{equation}\bigl[e^{-r_{ij}t}\bigr]\in \Jac S\simeq
\bigl(\oC^\ast\bigr)^g,\label{L^t}\end{equation} where the $r_{ij}$ are given in \eqref{aij}.

\subsection{Definite line bundles}

 If $F$ is real and $F(k-2)\not \in \Theta$, the space $H^0\bigl(S,F(k-1)\bigr)$ comes with a hermitian metric given by \eqref{form}.
 We compute this metric in the representation of sections given by \eqref{comps}.
Let $P_1,\dots,P_k$ and $R_1,\dots, R_k$ be two sets of polynomials defining sections $s,s^\prime$ of $F(k-1)$ via \eqref{comps}. Let us
write $\eta_i(\zeta)=z_i+2x_i\zeta-\bar{z}_i\zeta^2$. The formula \eqref{c0} gives (for any $\zeta$):
%%%%%%%%%%%%%%%%%%%%%%
\begin{equation}\langle s, s^\prime\rangle= (-1)^{k-1}\zeta^{k-1}\sum_{i=1}^k\frac{ P_{i}(\zeta)\overline{R_{i}\bigl(-1/\bar{\zeta}\bigr)}}{\prod_{j\neq
i} \bigl(\eta_i(\zeta)-\eta_j(\zeta)\bigr)}.\label{PR}\end{equation}
%%%%%%%%%%%%%%%%%%%%%%%%%%%%%%
Setting $\zeta=0$ we get
\begin{equation}\langle s, s^\prime\rangle= \sum_{i=1}^k\frac{
P_{i}(0)\overline{R_{i}(\infty)}}{\prod_{j\neq i} (z_i-z_j)},\label{lr}\end{equation} where $R_i(\infty)$ is the coefficient of
$\zeta^{k-1}$ in $R_i(\zeta)$.
\par
We can describe the set of definite (in the sense of definition \ref{pos}) line bundles of degree $g-1$. Let $\infty_i=S_i\cap
\pi^{-1}(\infty)$ and $0_i=S_i\cap \pi^{-1}(0)$, $i=1,\dots,k$. By assumption, these points are all distinct. Recall the definition
\eqref{Lambda} of the function $\Lambda_F$.
%%%%%%%%%%%%%%
\begin{proposition} The set of definite line bundles $F(k-2)$ on $S$ corresponds to the subset $\Jac^+ S$ of $\Jac S$ consisting of
equivalence classes (under the action \eqref{action}) $[\lambda_{ij}]\in \Jac S$ such that:
\begin{itemize}
\item[(1)] $\lambda_{ij}=\overline{\lambda_{ji}}$ for all $i\neq j$.
\item[(2)] For every $l=1,\dots,k$ the polynomial $Q_l(\zeta)=\Lambda_F(\infty,\dots,\infty,\zeta,0,\dots,0)$, where $\zeta$ occurs in
$l$-th place, has degree exactly $k-1$, i.e. its term of degree $k-1$ does not vanish.
\item[(3)] The sign of $\frac{Q_{l}(0)\overline{Q_{l}(\infty)}}{\prod_{j\neq l} (z_l-z_j)}$ is the same for each $l=1,\dots,k$.
\end{itemize}
\label{definite}\end{proposition}
\begin{proof} The first condition simply says that $F$ is real. The second condition with $l=k$ is equivalent, thanks to Proposition
\ref{zeros}, to $F(k-2)\not\in \Theta$. Condition (2) and Proposition \ref{zeros} also imply that, for each $l=1,\dots,k$, there is a
section $s_l$ of $F(k-1)$, unique up to a constant multiple, which vanishes at $\infty_1,\dots,\infty_{l-1},0_{l+1},\dots, 0_k$.
Proposition \ref{Ql} says that the restriction of $s_l$ to $S_l$ is the polynomial $Q_l$. Formula \eqref{c0} shows (as in the proof of
Proposition \ref{V_i}) that the $s_l$ are mutually orthogonal in the form \eqref{form} and that $\langle s_l, s_l \rangle$ is given by the
fraction in (3). Therefore \eqref{form} is definite if and only if all these fractions have the same sign.
\end{proof}

We do not know how many connected components $\Jac^+ S$ has. We think that there are at least two: one containing $L^t(k-2)$ for $t>0$ and the other one containing $L^t(k-2)$ for $t<0$. We have:
\begin{proposition} The form \eqref{form} is positive-definite on the connected component of $Jac^+ S$ containing $L^t(k-2)$ for $t>0$.\label{++}\end{proposition}
\begin{remark} From Corollary \ref{0}, we know that this component contains $F_t(k-2)$ for any real $F\in \Pic^0 S$ and sufficiently large $t$.\end{remark}
\begin{proof} Thanks to \eqref{L^t}, the bundle $L^t$ corresponds to $\bigl[e^{-r_{ij}t}\bigr]$ in $\Jac S$. Thus, it tends to the point $[0]$ on the boundary of $\Jac S$. We claim that there is a unitary basis $s^1(t),\dots,s^k(t)$ of  $L^t(k-1)$ which has a limit as $t\rightarrow +\infty$. This can be proved either via the asymptotic behaviour of solutions to Nahm's equations or by the argument used in the proof of Theorem \ref{bigdelta} below. From \eqref{match1}, we know that a limit section $(Q_1,\dots,Q_k)$ satisfies $Q_j(a_{ij})=0$ for any $i,j$, and hence it is generated by sections $s^l=(Q_1,\dots,Q_k)$ with $Q_l(\zeta)=\prod_{i\neq l}(\zeta-a_{il})$ and $Q_j\equiv 0$ for $j\neq l$. Thus, we shall be done as soon as we show that the norm of this $s^l$ is positive for each $l$. From the formula \eqref{PR} we need to compute $Q_{l}(\zeta)\overline{Q_{l}\bigl(-1/\bar{\zeta}\bigr)}$, which is:
$$ (-1)^{k-1}\zeta^{k-1}\prod_{i\neq l}(\zeta-a_{il})\bigl(-\frac{1}{\zeta}-\overline{a_{il}}\bigr)= \prod_{i\neq l}(\zeta-a_{il})\bigl(1+\zeta \overline{a_{il}}\bigr).$$
Since $a_{li}=-1/\overline{a_{il}}$, this can be rewritten as:
$$Q_{l}(\zeta)\overline{Q_{l}\bigl(-1/\bar{\zeta}\bigr)}=\prod_{i\neq l}\overline{a_{il}}\prod_{i\neq l}(\zeta-a_{il})\bigl(\zeta -{a_{li}})=\prod_{i\neq l}\frac{\overline{a_{il}}}{\bar{z}_i-\bar{z}_l} \prod_{i\neq
l} \bigl(\eta_l(\zeta)-\eta_i(\zeta)\bigr).$$
This, together with \eqref{PR} and \eqref{aij}, shows that $$\langle s^l, s^l \rangle= \prod_{i\neq l}\frac{x_i-x_l+r_{il}}{|z_i-z_l|^2},$$ which proves the result.
\end{proof}

\section{Hyperk\"ahler structure of adjoint orbits of $GL(k,\cx)$}

Let $z_1,\dots, z_k$ be distinct complex numbers and let us write $\tau=\diag(z_1,\dots,z_k)$. Let $O(\tau)$ be the adjoint $GL(k,\cx)$-orbit of $\tau$. It is a regular semi-simple orbit and Kronheimer \cite{Kr} shows that $O(\tau)$ admits a family of hyperk\"ahler structures, parameterised by $\oR^k$, such that the complex structure $I_1$ is the one of the complex adjoint orbit. We recall Kronheimer's construction.

\subsection{Orbits and Nahm's equations}
 The hyperk\"ahler structure on $O(\tau)$, given by a parameter  $(x_1,\dots,x_k)\in\oR^k$, is obtained by considering solutions to Nahm's equations \eqref{Nahm0} corresponding to the flow $L^t$ on the spectral curve \eqref{reducible}. The structure is seen better, if we  allow gauge freedom and introduce a fourth $\fU(k)$-valued function $T_0(t)$. Thus,
we consider the following variant of  Nahm's equations:
%%%
\begin{equation}\dot{T}_i+[T_0,T_i]+\frac{1}{2}\sum_{j,k=1,2,3}\epsilon_{ijk}[T_j,T_k]=0\;,\;\;\;\;i=1,2,3.\label{Nahm}\end{equation}
The functions $T_0,T_1,T_2,T_3$ are smooth, $\fU(k)$-valued and defined on the half-line $[0,+\infty)$. Moreover, if we write
\begin{equation} \tau=\tau_2+i\tau_3,\quad \tau_1=-i\diag(x_1,\dots,x_n),\label{tau}
\end{equation}
so that $\tau_1,\tau_2,\tau_3$ are diagonal skew-hermitian matrices, then we require that each $T_i(t)$, $i=1,2,3$, approaches $ \tau_i$ exponentially fast as $t\rightarrow +\infty$,  and that $T_0$ approaches $0$ also exponentially fast.
\par
The space of such solutions is acted upon by the gauge
group $\sG$ of $U(k)$-valued functions $g(t)$, with $g(0)=1$ and a diagonal limit, approached exponentially fast, at $+\infty$. The action is given by:
\begin{eqnarray} T_0&\mapsto & gT_0g^{-1}-\dot{g}g^{-1}\nonumber\\ T_i&\mapsto & gT_ig^{-1}\;,\;\;\qquad i=1,2,3.\label{Naction}\end{eqnarray}
%%%%%%%%%%
The hyperk\"ahler manifold $\sM_{\tau_1,\tau_2,\tau_3}$ is defined as the moduli space of solutions to \eqref{Nahm} satisfying these boundary
conditions  modulo the action of $\sG$.
\par
%%%
The tangent space at a solution $(T_0,T_1,T_2,T_3)$ can be identified with the space of solutions to the following system of linear
equations:
\begin{equation}\begin{array}{c} \dot{t}_0+[T_0,t_0]+[T_1,t_1]+[T_2,t_2]+[T_3,t_3]=0,\\
\dot{t}_1+[T_0,t_1]-[T_1,t_0]+[T_2,t_3]-[T_3,t_2]=0,\\
\dot{t}_2+[T_0,t_2]-[T_1,t_3]-[T_2,t_0]+[T_3,t_1]=0,\\
\dot{t}_3+[T_0,t_3]+[T_1,t_2]-[T_2,t_1]-[T_3,t_0]=0.\end{array}\label{tangent}\end{equation} The first equation is the condition that
$(t_0,t_1,t_2,t_3)$ is orthogonal to the infinitesimal gauge transformations and the remaining three are linearisations of \eqref{Nahm}.
\par
$\sM_{\tau_1,\tau_2,\tau_3}$ carries a hyperk\"ahler metric defined by
\begin{equation}\|(t_0,t_1,t_2,t_3)\|^2=-\sum_{i=0}^3\int_{0}^{+\infty}\tr t_i^2(s) ds \label{metric},\end{equation}
and the three anticommuting structures $I_1,I_2,I_3$ are given by the left multiplication of $t_0+t_1i+t_2j+t_3k$ by $i,j,k$ ($i,j,k$ - the standard basis of imaginary quaternions).
\par
Kronheimer shows that the map
$$\bigl(\sM_{\tau_1,\tau_2,\tau_3},I_1\bigr)\rightarrow O(\tau_2+i\tau_3),\quad \bigl(T_0(t),T_1(t),T_2(t),T_3(t)\bigr)\mapsto T_2(0)+iT_3(0)$$
is a biholomorphism.  In fact, if we identify, as usual, the complex structures with $\oP^1$, then the map
$$\phi_\zeta:\bigl(\sM_{\tau_1,\tau_2,\tau_3},I_\zeta\bigr)\rightarrow O(\tau_\zeta), \quad \bigl(T_0(t),T_1(t),T_2(t),T_3(t)\bigr)\mapsto A(0,\zeta),$$
$$ \tau_\zeta=(\tau_2+i\tau_3)+2i\tau_1\zeta +(\tau_2-i\tau_3)\zeta^2,$$ $$ A(0,\zeta)=(T_2(0)+iT_3(0))+2iT_1(0)\zeta +(T_2(0)-iT_3(0))\zeta^2,$$
 a biholomorphism as long as $\tau_\zeta$ is a regular matrix. For other complex structures, $\sM_{\tau_1,\tau_2,\tau_3}$ is a vector bundle over a flag manifold \cite{Biq}.

\subsection{Twistor lines via $\Jac S$} We wish to make a connection between the hyperk\"ahler geometry of $\sM_{\tau_1,\tau_2,\tau_3}$ and the algebraic geometry of a reducible curve $S$. The curve $S$ is given by the equation
$$\det\left(\eta\cdot 1-(\tau_2+i\tau_3)-2i\tau_1\zeta -(\tau_2-i\tau_3)\zeta^2\right)=0.$$
It is of the form \eqref{reducible} with $z_i$ (resp. $x_i$) being the diagonal entries of $\tau_2+i\tau_3$ (resp. $i\tau_1$). In order to make use of previous sections, we assume that $S$ has only nodes, i.e. no three of the points $(x_i, \re z_i, \im z_i)\in \oR^3$, $i=1,\dots,k$, are collinear.

The action of $\sG$ allows one to eliminate the component $T_0$ in the definition of $\sM_{\tau_1,\tau_2,\tau_3}$. Thus, as a manifold, $\sM_{\tau_1,\tau_2,\tau_3}$ is diffeomorphic to triples $(T_1,T_2,T_3)\in \u(n)\otimes \oR^3$ such that the solution $\bigl(T_1(t),T_2(t),T_3(t)\bigr)$ of the Nahm's equations \eqref{Nahm0} with the initial condition $T_i(t)=T_i$, $i=1,2,3$, exists for all $t\geq 0$ and $\lim_{t\rightarrow +\infty} T_i(t)$ exists and is conjugate to $\tau_i$, $i=1,2,3$. We know from section \ref{lb} that, as long as $A(\zeta)=(T_2+iT_3)+2iT_1\zeta +(T_2-iT_3)\zeta^2$ is a regular matrix for every $\zeta$, the $U(n)$-conjugacy class of $(T_1,T_2,T_3)$ is identified with an element of $\Jac^+ S$. We do not obtain all points of $\Jac^+ S$ this way: only those $F$ such that $F_t(k-2)\not\in \Theta$ for all $t\geq 0$. We shall denote this subset by $B$. We expect that $B$ is a connected component of $\Jac^+ S$.
\par
We remark that the remaining (i.e. non-regular) triples $(T_1,T_2,T_3)\in \u(n)\otimes \oR^3$ correspond to the closure $\overline{B}$ of $B$ in $\overline{\Jac_\oR S} -\overline{\Theta}_\oR$, where $\overline{\Jac S}$ is the compactified Jacobian and $\overline{\Theta}$ is the closure of $\Theta$ in $\overline{\Jac S}$ (cf. Remarks \ref{comp} and \ref{compactified}).
\par
It is clear from the above discussion that $\overline{B}$ is identified with $\sM_{\tau_1,\tau_2,\tau_3}/U(k)$ and, hence, with $O(\tau)/U(k)$, where $\tau=\tau_2+i\tau_3$. We wish to describe the fibration $\sM_{\tau_1,\tau_2,\tau_3}\rightarrow \overline{B}$ in greater detail.
\par
As in section \ref{lb}, we denote by $V$ a vector bundle over $B$, whose fibre at $F(k-2)$ is $H^0\bigl(S,F(k-1)\bigr)$. Proposition \ref{++} implies that the metric \eqref{form} is  positive-definite on $V$. Let us denote by
 $ U(V)$ the corresponding bundle of unitary frames of $V$.  Topologically, $U(V)$ is a trivial bundle (we can find a section as in the proof of Proposition \ref{definite}). Let $PU(V)$ be the quotient of $U(V)$ by the centre of $U(k)$.  We can restate Beauville's theorem \ref{Beauville} in the following form:
%%%%%%%%%%%%%%%%%%%%%%%%
\begin{proposition} There is a canonical $U(k)$-equivariant diffeomorphism between the restriction of $PU(V)$ to $B$ and the corresponding open dense subset of $\sM_{\tau_1,\tau_2,\tau_3}$.\end{proposition}
\begin{proof} Let $F(k-2)\in B$. If $\psi^1,\dots,\psi^k$ is a unitary basis of $H^0\bigl(S,F(k-1)\bigr)$, then we obtain a a triple $(T_1,T_2,T_3)$ in $\sM_{\tau_1,\tau_2,\tau_3}$ via \eqref{conjugate}. This map is smooth and $U(k)$-equivariant. The induced map from $PU(V)$ is injective by the argument used in the proof of Proposition \ref{realbundles}. The inverse mapping is provided by \eqref{bundle}.
\end{proof}

Thus, Kronheimer's hyperk\"ahler structure lives naturally on $PU(V)$. We are going to describe the twistor lines from this point of view.
\par
First let us fix notation. Let $Q=[Q_{ij}]$ be a matrix of complex polynomials of degree $n$. We write $\sigma(Q)$ for the matrix of polynomials defined by $\sigma(Q)_{ij}(\zeta)=(-1)^n\overline{Q_{ji}(-1/\bar{\zeta})}$.

\begin{theorem} The bundle $U(V)$ is isomorphic to an open $U(k)$-invariant subset $\sV(S)$ of the set of all $k\times k$ matrices $Q$ of polynomials of degree $k-1$ which satisfy:
\begin{equation}\sigma(Q)DQ=1,\label{QDQ}\end{equation} where
$$D=\diag\left(\frac{1}{\prod_{j\neq
i} \bigl(\eta_i(\zeta)-\eta_j(\zeta)\bigr)}\right).$$
The action of $U(k)$ is given by right multiplication on elements of $Q$ and the projection $\sV(S)\rightarrow B$ is $Q\mapsto [\lambda_{mn}]$, where $\lambda_{mn}=\frac{Q_{ns}}{Q_{ms}}(a_{mn})$ for any $s$. \label{U(V)}\end{theorem}
\begin{proof} The map from $U(V)$ to $\sV(S)$ is clear: represent a unitary basis of $H^0\bigl(S,F(k-1)\bigr)$ in the form \eqref{comps}. This gives a quadratic matrix of polynomials satisfying \eqref{QDQ}, thanks to \eqref{PR}.
\par
Conversely, suppose that $Q$ is a matricial polynomial satisfying \eqref{QDQ}.
We need the following lemma:
\begin{lemma} Let $Z_1(\zeta),\dots, Z_k(\zeta)$ be polynomials of degree $2k-2$ such that
\begin{equation}\sum_{s=1}^k \frac{Z_s(\zeta)}{\prod_{l\neq
s} \bigl(\eta_s(\zeta)-\eta_l(\zeta)\bigr)}=\text{const}.\label{Zcon}\end{equation}
Then, there exists a section $s$ of $\sO_S(2k-2)$  such that $Z_i(\zeta)$ is a restriction of $s$ to $S_i$.\end{lemma}
  Indeed, suppose that $Z=\bigl(Z_1(\zeta),\dots, Z_k(\zeta)\bigr)$ satisfies this condition. If the constant is nonzero, we can replace $Z_i$ by $Z_i-c\eta_i^{k-1}$, where $c$ is the constant in question. Thus, we can assume that the $Z_i$ satisfy the homogeneous version of this equation.
It is clear that any section of $\pi^\ast\sO(2k-2)$, which must be of the form \eqref{c_0}, satisfies the above condition. Therefore we can add to $Z$ a section of the form \eqref{c_0} in order to assume that $Z_i(a_{ij})$  are all non-zero. Equation \eqref{Zcon} implies, in particular, that the residue at $a_{ij}$ of the left-hand side vanishes. The only terms which have a pole at $a_{ij}$ are $s=i$ and $s=j$. We observe that, if we replace all $Z_i$ in \eqref{Zcon} by $1$, we obtain $0$ (this can be viewed as the product of the last row of $V^{-1}$  with the first column of $V$, where $V$ is the Vandermonde determinant in $\eta_1(\zeta),\dots,\eta_k(\zeta)$). Therefore $Z_i(a_{ij})=Z_j(a_{ij})$, which proves the lemma.
\par
We return to the proof of the proposition.
We consider the open subset $\sV(S)$, where the $i$-th column of $Q$, $i=1,\dots,k$, defines a section of $F_i(k-1)$ for some $F_i\in \Pic^0 S$ with $F_i(k-2)\not\in \Theta$. The $i$-th row $R_{i1},\dots, R_{ik}$ of $\sigma(Q)$ is then, because of \eqref{sigma0}, a section of $\sigma(F_i)(k-1)$. Thus,  $R_{i1}Q_{1j},\dots, R_{ik}Q_{kj}$ represents a section of $F_j\otimes \sigma( F_i)(2k-2)$, and, from the above lemma,  $F_j\otimes \sigma( F_i)\simeq \sO$ for all $i,j$. Taking first $i=j$  and then all $i,j$ shows that all $F_i$ are isomorphic to one bundle $F$ which is real. The matrix $Q$ becomes now, thanks to \eqref{QDQ}, a unitary basis of sections of $F(k-1)$ which proves that $F\in \Jac^+ S$ and the metric is positive-definite.
\end{proof}

 It is clear how the hyperk\"ahler structure of $\sM_{\tau_1,\tau_2,\tau_3}$ looks on  $\sV(S)/U(1)$. For every $\zeta_0$, which is not one of the $a_{mn}$, we have a map from $\sV(S)$ to the adjoint orbit of $ \tau_{\zeta_0}=(\tau_2+i\tau_3)+2i\tau_1\zeta_0 +(\tau_2-i\tau_3)\zeta_0^2$, given by
$$ Q\longmapsto Q(\zeta_0)^{-1}\tau_{\zeta_0} Q(\zeta_0)$$
(we observe that \eqref{QDQ}  guarantees that $Q(\zeta)$ is invertible if $\zeta$ is not one of the $a_{mn}$).
The pull-backs of complex structures of $O(\tau_{\zeta_0})$, for this dense subset of $\zeta_0$, to $\sV(S)/U(1)\simeq PU(V)$  generate the hypercomplex structure of $\sV(S)/U(1)$ (isomorphic to that of $\sM_{\tau_1,\tau_2,\tau_3}$).

\begin{remark} To extend this description to all of $\sM_{\tau_1,\tau_2,\tau_3}$, we have to consider a bigger subset of polynomial matrices satisfying \eqref{QDQ}: those, whose columns fibre over a point in $\overline{\Jac S}-\overline{\Theta}$.\end{remark}

\section{K\"ahler potentials}

\subsection{K\"ahler potentials in terms of Nahm's data}
 It is well-known (see, e.g. \cite{Hit2}) that the function
 $$K= -\frac{1}{2}\int_0^\infty \tr \bigl(T_1^2(t)+T_2^2(t)\bigr)dt$$
 is a $U(n)$-invariant K\"ahler potential for $\bigl(\sM_{0,0,\tau_3},I_1\bigr)$, i.e. $-i\partial_{I_1}\bar{\partial}_{I_1}K$ is the K\"ahler form $\omega_1=g(I_1\cdot,\cdot)$ of the hyperk\"ahler metric $g$ given by \eqref{metric}. The reducible spectral curve $S$ corresponding to $\tau_1=\tau_2=0$ has, however, points of multiplicity greater than two, and, so, this $S$ does not fit in with the results of the previous sections. On the other hand, the cohomology class $[\omega_1]$ is identified with $\tau_1$ \cite{Kr}, and, hence,  $\bigl(\sM_{0,\tau_2,\tau_3},I_1,\omega_1\bigr)$, being Stein, must admit a $U(k)$-invariant K\"ahler potential. We have:
 %%%%
 \begin{proposition} The function
 \begin{equation} K= -\frac{1}{2}\int_0^\infty \tr \bigl(T_1^2(t)+T_2^2(t)-\tau_2^2\bigr)dt\label{K}\end{equation}
 is a K\"ahler potential for $\bigl(\sM_{0,\tau_2,\tau_3},I_1,\omega_1\bigr)$.\label{KK}\end{proposition}
 \begin{proof} This can be proved via a direct calculation, but it is more instructive to use the fact that $\sM_{\tau_1,\tau_2,\tau_3}$ is the (infinite-dimensional) hyperk\"ahler quotient\footnote{Alternatively, the same argument can be applied to $\sM_{\tau_1,\tau_2,\tau_3}$ as a finite-dimensional hyperk\"ahler quotient of the manifold $M_{U(k)}(c)$, with its K\"ahler potential, described in section 4 of \cite{RB}} of the space $\sA$ of all quadruples $\bigl(T_0(t),T_1(t),T_2(t),T_3(t)\bigr)$ satisfying the prescribed boundary conditions at $+\infty$, with its flat hyperk\"ahler structure, by the gauge group $\sG$. The moment map equations  are the three equations \eqref{Nahm}. The function $K$ is finite, since $\tau_1=0$, and is a K\"ahler potential for $\bigl(\sA,I_1,\omega_1\bigr)$ by direct computation. It is still a K\"ahler potential on the complex submanifold $\sB$ defined as the space of solutions to \eqref{Nahm} with $i=2$ and $i=3$. $\bigl(\sM_{0,\tau_2,\tau_3},I_1,\omega_1\bigr)$ is the K\"ahler quotient of $\bigl(\sB,I_1,\omega_1\bigr)$ by $\sG$ and we shall use the following lemma:
 \begin{lemma} Let $K$ be a globally defined $G$-invariant K\"ahler potential on a K\"ahler manifold $(M,I,\omega)$ on which the action of $G$ is Hamiltonian and holomorphic. Let $\mu:M\rightarrow \fG^\ast$ be a moment map on $M$ defined by $\langle\mu,\rho\rangle =
dK(I\hat{\rho})$, where $\hat{\rho}$ is the fundamental vector field induced by $\rho\in\fG$.
 \par
  Then $K$ descends to a K\"ahler potential for the K\"ahler quotient  $\mu^{-1}(0)/G$.\end{lemma}
To see that this proves the proposition, we have to check that $dK$, where $K$ is given by \eqref{K}, evaluated on vectors $I_1\hat{\rho}$, $\rho\in\Lie \sG$, gives the moment map with respect to $\omega_1$, i.e. \eqref{Nahm} holds for $i=1$. We compute $I_1\hat{\rho}$ by differentiating \eqref{Naction}:
$$I_1\hat{\rho}= \bigl( -[\rho,T_1], [\rho,T_0]-\dot{\rho}, -[\rho,T_3], [\rho,T_2]\bigr).$$
Then
$$dK(I_1\hat{\rho})=-\int_0^\infty \tr \bigl(T_1([\rho,T_0]-\dot{\rho})-T_2[\rho,T_3]\bigr)=-\int_0^\infty \tr\rho\bigl(\dot{T}_1+[T_0,T_1]+[T_2,T_3]\bigr),$$
 where we have used the fact that $\tau_1=0$. It remains to prove the lemma.
\newline
  {\em Proof of Lemma}. Let $\pi:\mu^{-1}(0)\rightarrow \mu^{-1}(0)/G$ be the projection. Let $\bar{I}$ and $\bar{K}$ denote the induced complex structure and induced map $K$ on $\mu^{-1}(0)/G$. In other words $K=\bar{K}\circ\pi$ and $\bar{I}v$ is $\pi(Iv^h)$, where $v^h$ denotes the horizontal lift of $v$. We compute:
  $$ \bar{I}d\bar{K}(v)=-dK(Iv^h)=IdK(v^h).$$
  By assumption, $IdK(\hat{\rho})=0$ at points of $\mu^{-1}(0)$ and, hence,
  $IdK(v^h)=IdK(\tilde{v})$ for any other lift $\tilde{v}$ of $v$. Thus, $ \bar{I}d\bar{K}(\pi(u))=IdK(u)$ and taking the exterior derivative of both sides proves the lemma.
\end{proof}

 We can identify the K\"ahler potential \eqref{K} for the hyperk\"ahler metric on the orbit $O(\tau)$ using ideas of Hitchin \cite{Hit2}. We need to extend Theorem 2 on p. 56 in \cite{Hit2} to fully reducible curves.

\subsection{Hitchin's theorem for reducible curves}

For a matricial polynomial $A(\zeta)=A_0+A_1\zeta+A_2\zeta^2$, representing a line bundle in $J^{g-1}(S)-\Theta$, the expression
\begin{equation} \Delta=\tr\left( A_0A_2-\frac{1}{4} A_1^2\right)\label{Delta} \end{equation}
is $GL(k,\cx)$ invariant and, hence, it defines a holomorphic function on $J^{g-1}(S)-\Theta$. Hitchin \cite{Hit2} shows that $\Delta$
extends to a meromorphic function on $J^{g-1}(S)$ and, {\em for a smooth $S$}, relates it to the theta function of $S$. Our next aim is to
extend this to fully reducible spectral curves.
\par
Let $S$ be a curve given by \eqref{reducible}. Recall from section \ref{spec} that the line bundle $L^t$, which on $T\oP^1$ has the
transition function $e^{-t\eta/\zeta}$, induces the element $\bigl[e^{-r_{ij}t}\bigr]$ in $\Jac S\simeq \bigl(\cx^\ast\bigr)^g$. This
induces a flow $[\lambda_{ij}]\rightarrow \bigl[e^{-r_{ij}t}\lambda_{ij}\bigr]$ and, hence, a vector field $X$ on $\Jac S$. For a function $f$ on $\Jac S$ we write
$$\frac{df}{dt}=Xf=\frac{d}{dt}f\bigl[e^{-r_{ij}t}(\cdot)\bigr].$$
\par
Let us also reformulate $\Delta$. We can write $A_0=T_2+iT_3$, $A_1=2iT_1$, $A_2=T_2-iT_3$ for complex matrices $T_1,T_2,T_3$. Then
$$\Delta = \tr \left(T_1^2+T_2^2+T_3^2\right).$$
We recall that $\tau_1,\tau_2,\tau_3$ are diagonal skew-hermitian matrices defined by $\tau_2+i\tau_3=\diag(z_1,\dots,z_k)$, $i\tau_1=\diag(x_1,\dots,x_k)$.
\par
 We now have the following analogue of Theorem 2 in \cite{Hit2}:
%%%%%%%%%%%%%%%%%
\begin{theorem} Let $S$ be a reducible spectral curve given by \eqref{reducible} with nodes only. Then, the following equality holds on $\Jac S$:
\begin{equation} \Delta-\tr\left(\tau_1^2+\tau_2^2+\tau_3^2\right)=\frac{3}{2}\frac{d^2}{dt^2}\log \vartheta_{1,0},\label{d^2}\end{equation}
where $\Delta$ is given by \eqref{Delta} and $\vartheta_{1,0}$ is defined in Proposition \ref{p,q}.\label{bigdelta}\end{theorem}
\begin{proof}
We follow Hitchin's strategy with suitable modifications due to non-compactness of $\Jac S$. The computation of polar parts of $\Delta$ and
$\vartheta_{1,0}$, given in \cite{Hit2}, remains valid and, hence, we know that $\Delta-\frac{3}{2}\frac{d^2}{dt^2}\log \vartheta_{1,0}$ is
a holomorphic function on $\Jac S$. To prove that it is equal to $\tr\left(\tau_1^2+\tau_2^2+\tau_3^2\right)$, we need to study its
asymptotic behaviour.
\par
Let $[\lambda_{ij}]\in \Jac S$ and consider its orbit under the flow $\bigl[e^{-r_{ij}t}\lambda_{ij}\bigr]$ ($t\in\cx$). We assume at first that all $r_{ij}$ are rational, so that this orbit is a closed complex submanifold of $\Jac S$
isomorphic to $\cx^\ast$. We claim that both $\Delta-\tr\left(\tau_1^2+\tau_2^2+\tau_3^2\right)$ and $\frac{d^2}{dt^2}\log
\vartheta_{1,0}$, restricted to the flow orbit, approach zero as $\re t\rightarrow \pm \infty$. This, of course, will prove the result for rational $r_{ij}$ and, hence, by continuity, for all $S$. We deal first with
$\vartheta_{1,0}$. We have:
$$\frac{d^2}{dt^2}\log
\vartheta_{1,0}=\frac{1}{\vartheta}\frac{d^2\vartheta}{dt^2}-\left(\frac{1}{\vartheta}\frac{d\vartheta}{dt}\right)^2,$$
%%%%%%%%%%%%%%%%
We now observe that $\vartheta_{1,0}$ restricted to the flow orbit approaches $\vartheta_{1,0}([0])$ exponentially fast as $\re t
\rightarrow +\infty$. Moreover, we know from Corollary \ref{0}  that $\vartheta_{1,0}([0])\neq 0$ and, hence, $\frac{d^2}{dt^2}\log
\vartheta_{1,0}$ approaches $0$ exponentially fast, as $\re t \rightarrow +\infty$. To deal with $t\rightarrow -\infty$, we observe that
dividing $\vartheta_{1,0}$ by $\prod_{i\neq j}\lambda_{ij} e^{-r_{ij}t}$ yields $\vartheta_{0,1}$ restricted to the orbit
$\bigl[e^{r_{ij}t}\lambda_{ij}^{-1}\bigr]$ but it does not change $\frac{d^2}{dt^2}\log \vartheta_{1,0}$. Thus, the argument used for
$\vartheta_{1,0}$ can be now applied to $\vartheta_{0,1}$ and, so, $\frac{d^2}{dt^2}\log \vartheta_{1,0}$ approaches zero exponentially
fast, as $\re t \rightarrow \pm\infty$.
\par
We now show that $\Delta$ approaches $\tr\left(\tau_1^2+\tau_2^2+\tau_3^2\right)$ as $\re t\rightarrow \pm \infty$. Let $F\in \Pic^0 S$ be
the line bundle represented by $[\lambda_{ij}]$. From section \ref{lb} we know that there is an anti-holomorphic isomorphism
$$H^0\bigl(S,F(k-1)\bigr)\longrightarrow H^0\bigl(S, \sigma( F)(k-1)\bigr)$$
%%%%%%%%%%%%%%%
 and, hence, $\Delta - \tr\left(\tau_1^2+\tau_2^2+\tau_3^2\right)$
 approaches $0$ as $\re t\rightarrow -\infty$ on the orbit $F_t(k-1)$ if and only if it  approaches $0$ as $\re t\rightarrow +\infty$ on the
orbit  $\sigma( F)(k-1)$. Therefore we only need to deal with $\re t\rightarrow +\infty$. Let $t_n$ be a sequence of complex
numbers with $\re t_{n}\rightarrow +\infty$. The fact that $\vartheta_{1,0}([0])\neq 0$ implies that $F_{t_n}(k-2)\not \in \Theta$ for $n$
large enough. Therefore $F_{t_n}(k-1)[-D_0]\not\in \Theta$ for $n$ large enough, where $D_0=0_1+\dots+0_k$ is the divisor $\pi^{-1}(0)$
(with $0_i\in S_i$), and we can find a basis of sections $s^1_n,\dots,s^k_n$ of $F_{t_n}(k-1)$ such that $s^l_n$ vanishes at all $0_i$,
$i\neq l$. We represent this basis in the form \eqref{comps}, so that $s^l_n$ is given by polynomials $Q_{1l}^n(\zeta),\dots,
Q_{kl}^n(\zeta)$ satisfying the matching conditions $Q_{jl}^n(a_{ij)}=e^{-r_{ij}t_n}\lambda_{ij}Q_{il}^n(a_{ij})$ for all $i\neq j$ and all
$l,n$. In addition, $Q_{il}^n(0)=0$ for $i\neq l$ and all $n$. Let now $R$ be such that all $a_{ij}$ are contained in the disc $B_R\subset
\oC$ of radius centred at $0$ and of radius $R$. Define $M_l^n$ as the sup norm of $\{Q_{1l}^n(\zeta),\dots,Q_{kl}^n(\zeta)\}$ over $B_R$.
For a fixed $l$, there a subsequence of $n$ such that $Q_{jl}^n(\zeta)/M_l^n$ converges for every $j=1,\dots,k$, and at least one of these
has a non-zero limit. Taking the limit of the matching conditions satisfied by the $Q_{jl}^n$ shows that their limits $Q_{jl}$ satisfy
$Q_{jl}(a_{ij})=0$ for every $i\neq l$. Since these are polynomials of degree at most $k-1$ and, in addition, $Q_{jl}(0)=0$ for $j\neq l$,
we have that $Q_{jl}\equiv 0$ for $j\neq l$ and $Q_{ll}(\zeta)=C_l\prod_{i\neq l} (\zeta-a_{il})$ for a nonzero $C_l$. We can find a
subsequence such that these limits exist for every $l=1,\dots,k$. Now recall from \eqref{conjugate} that a matricial polynomial
 corresponding to $F_{t_n}(k-1)$ is obtained by conjugating a diagonal matrix by the matrix $\bigl[Q_{ij}^n(\zeta)\bigr]$. Another
 matricial polynomial $A_n(\zeta)$ (conjugate to the previous one) is obtained by conjugating by $\bigl[Q_{ij}^n(\zeta)/M_j^n\bigr]$. We showed
above that this latter matrix has a diagonal limit (with non-vanishing diagonal entries). Therefore a subsequence of $A_n(\zeta)$ has a diagonal limit equal
to $(\tau_2+i\tau_3)+2i\tau_1\zeta+(\tau_2-i\tau_3)\zeta^2$. This finishes the proof.
\end{proof}

\begin{corollary} For $i=1,2,3$, $\tr T_i^2-\tr \tau_i^2=\frac{1}{2}\frac{d^2}{dt^2}\log \vartheta_{1,0}$.\end{corollary}
\begin{proof} The functions $\tr T_i^2 -T_j^2$, $i,j=1,2,3$, all arise from $\tr A(\zeta)^2$, which is constant on the Jacobian. \end{proof}

\subsection{A formula for the K\"ahler potential of Kronheimer's metric}

We can now identify the K\"ahler potential \eqref{K} for the complex structure $I_1$ of the hyperk\"ahler manifold $\sM_{0,\tau_2,\tau_3}$. As a complex manifold $\bigl(\sM_{0,\tau_2,\tau_3},I_1\bigr)$ is the adjoint orbit $O(\tau_2+i\tau_3)$ of $GL(k,\cx)$. We assume that the diagonal entries $z_1,\dots,z_k$ of $\tau_2+i\tau_3$ are all distinct, so that the orbit is regular semisimple. In addition, we require that the spectral curve $S$, given by \eqref{reducible} with $x_1=\dots=x_k=0$, has only ordinary double points, i.e. we assume that no three points $z_i\in \cx$ are collinear. With these assumptions, we know that $K$, being $U(k)$-invariant, is given by a function on the subset $B$ of $\Jac_\oR S-\Theta_\oR$. We have:
%%%%%%%%%%%%%%%%%%%%%%
\begin{theorem} The K\"ahler potential $K$, defined by \eqref{K}, of $\bigl(\sM_{0,\tau_2,\tau_3},I_1,\omega_1\bigr)$ is given by
$$ K(T_0,T_1,T_2,T_3)= \frac{1}{2}\frac{\dot{\vartheta}_{1,0}(\lambda)}{\vartheta_{1,0}(\lambda)},$$
where $\lambda$ is the point of $\Jac S$ corresponding to the triple $\bigl(T_1(0),T_2(0),T_3(0)\bigr)$ and $\vartheta_{1,0}$ is defined in Proposition \ref{p,q}.\end{theorem}
\begin{proof} We just follow \eqref{K}, i.e.  integrate,  using the last corollary and the fact, already observed in the proof of Theorem \ref{bigdelta}, that $\dot{\vartheta}_{1,0}/\vartheta_{1,0}$ decays exponentially fast as $t\rightarrow +\infty$.\end{proof}

\begin{example} We compute the K\"ahler potential \eqref{K} in the case when $k=2$ and $z_2=-z_1=R/2$, i.e. we give yet another computation for the Eguchi-Hanson metric (cf. \cite{Hit2}). The theta function \eqref{theta} is $\lambda_{12}\lambda_{21}- \mu_{12}\mu_{21}$ and, after setting $\mu_{12}=\mu_{21}=1$ and quotienting by \eqref{action}, we have $\vartheta_{1,0}(\gamma)=\gamma-1$, where $\gamma= \lambda_{12}\lambda_{21}$ is the coordinate on $\cx^\ast\simeq \Jac S$. Thus, $\vartheta_{1,0}(\gamma,t)=e^{-2Rt}\gamma-1$ and $\ddot{\vartheta}=-2R\dot{\vartheta}$.
This means that $\frac{d^2}{dt^2}\log
\vartheta_{1,0}=-2R\frac{d}{dt}\log
\vartheta_{1,0}-\bigl(\frac{d}{dt}\log
\vartheta_{1,0}\bigr)^2$ and, writing $f=\frac{d}{dt}\log
\vartheta_{1,0}$ and $X=T_2(0)+iT_3(0)$, we get from the last corollary: $f^2+2Rf-\tr XX^\ast+R^2/2=0$. Therefore $f=-R+\sqrt{R^2/2+\tr XX^\ast}$ (as $f$ must be nonnegative on the orbit). The potential $K$ is half of $f$.
\end{example}

\end{document}